\documentclass[11pt,reqno]{amsart}
\usepackage[DIV=12, oneside]{typearea}
\usepackage[utf8]{inputenc}
\usepackage[T1]{fontenc}
\usepackage[english]{babel}
\usepackage{amsmath}
\usepackage{esint}
\usepackage{amssymb, amsthm, bbm, dsfont}
\usepackage{enumerate}
\usepackage{color, mathrsfs}
\usepackage{todonotes}
\usepackage{thmtools}

\usepackage{latexsym}
\usepackage{epsfig,mathrsfs}
\usepackage{hyperref}
\hypersetup{
    colorlinks=true,       
    linkcolor=blue,          
    citecolor=blue,        
    filecolor=blue,      
    urlcolor=cyan  
}

\usepackage{array}
\usepackage{booktabs}
\usepackage{stmaryrd}

\declaretheorem[name=Theorem, parent=section]{Thm}
\declaretheorem[name=Lemma, sibling=Thm]{Lem}
\declaretheorem[name=Proposition, sibling=Thm]{Prop}

\theoremstyle{definition}
\declaretheorem[name=Definition, sibling=Thm]{Def}

\newcommand{\R}{{\ensuremath{\mathbb{R}}}}
\newcommand{\N}{{\ensuremath{\mathbb{N}}}}

\renewcommand{\P}{\ensuremath{\mathbb{P}}}
\renewcommand{\dj}{d\kern-0.4em\char"16\kern-0.1em}
\newcommand{\E}{\ensuremath{\mathbb{E}}}

\newcommand{\cL}{  \mathcal{L}}

\newcommand{\cS}{  \mathcal{S}}

\renewcommand{\epsilon}{\varepsilon}

\renewcommand{\d}{\, \textnormal{d}}

\addto\extrasenglish{%

}

\makeatletter 
\def\Links{\tagsleft@true}\def\Rechts{\tagsleft@false} 
\makeatother 

\begin{document}
\numberwithin{equation}{section}
\bibliographystyle{amsalpha}

\title{Intrinsic scaling properties for
nonlocal operators}
\begin{abstract}
We study growth lemmas and questions of regularity for generators of Markov
processes. The generators are allowed to have an arbitrary order of
differentiability less than $2$. In general, this order is
represented by a function and not by a number. The approach enables a
careful study of regularity issues up to the phase boundary between
integro-differential (positive order of differentiability) and integral
operators (nonnegative order of differentiability). The proof is based
on intrinsic scaling properties of the underlying operators and stochastic
processes. 
\end{abstract}

\author{Moritz Kassmann}
\author{Ante Mimica}

\address{Fakultät für Mathematik\\
Universität Bielefeld \\
Postfach 100131 \\
D-33501 Bielefeld \\
Germany}
\email{moritz.kassmann@uni-bielefeld.de}

\address{Department of Mathematics \\ University of Zagreb \\ Bijeni\v{c}ka
cesta 30 \\ 10000 Zagreb \\ Croatia}
\email{amimica@math.hr}

\thanks{Research supported by German Science Foundation (DFG) via SFB
701. Research supported by MZOS grant 037-0372790-2801.}



\maketitle

\section{Introduction}

One key argument in the regularity theory of differential equations of second
order is the so called growth lemma. Here is an example which is by now
classical. Let $A$ be an elliptic operator of second
order, e.g. $Au = \sum_{i,j} a_{ij}(\cdot)
\tfrac{\partial}{\partial x_{i}} \tfrac{\partial}{\partial x_{j}} u$ for $u:\R^d
\to \R$ where
$(a_{ij}(\cdot))_{i,j}$ is uniformly positive definite and bounded. One
could
also consider nonlinear examples. The following growth lemma holds true
in many cases:

\begin{Lem}\label{lem:growth} There is a constant $\theta \in (0,1)$ such that,
if $R>0$ and $u:\R^d \to \R$ with
\begin{align*}
 &-A u \leq 0 \text{ in } B_{2R} \,, \qquad u \leq 1 \text{ in }
B_{2R}, \qquad |(B_{2R} \!\setminus\! B_R) \cap \{u \leq 0\}| \geq \tfrac12
|B_{2R} \!\setminus\! B_R|\,,
\end{align*}
then $u \leq 1-\theta$ in $B_R$.
\end{Lem}

Such lemmas are systematically studied and applied in \cite{Lan71}. Their
importance is underlined in the article \cite{KrSa79}, in
which the authors establish a priori bounds for elliptic
equations of second order with bounded measurable coefficents. Nowadays they
form a standard
tool for the study of various questions of nonlinear partial differential
equations of second order, cf. \cite{CaCa95} and \cite{DBGV12}. Note that the
property
formulated in \autoref{lem:growth} is also referred to as expansion
of positivity which describes the corresponding property for $1-u$.

In the case of a linear differential operator $A$ the above
lemma can be established with the help of the Markov process it generates. Let
$X$ be the strong Markov process associated with the operator $A$, i.e. we
assume that the martingale problem has a unique solution. Denote by $T_A,
\tau_A$ the hitting resp. exit time for a measurable set $A \subset \R^d$
and by
$\P_x$ the measure on the path space with $\P_x(X_0=x)=1$. The following
property is a key to the above growth lemma.

\begin{Prop}\label{prop:hitting-prop}
There is a constant $c \in (0,1)$ such that for every $R>0$ and every
measurable set $A \subset B_{2R} \!\setminus\!B_R$ with $|(B_{2R} \!\setminus\!
B_R) \cap A| \geq \tfrac12 |B_{2R} \!\setminus\! B_R|$ and $x \in B_R$
\begin{align}\label{eq:hitting-estim}
\P_x (T_A < \tau_{B_{2R}}) \geq c \,.
\end{align}
\end{Prop}

The aim of this work is to establish a result like \autoref{prop:hitting-prop}
and regularity estimates for a general class of operators and stochastic
processes. The article
\cite{KrSa79} deals with a very specific case: operators of second order.
Another very specific case, operators of fractional order $\alpha \in
(0,2)$, is treated in \cite{BL1}. Therein it is shown that
\autoref{prop:hitting-prop} holds true for jump
processes $X$ generated by integral
operators ${\cL}\colon C^2_b(\R^d)\rightarrow C(\R^d)$ of the form
\begin{align}\label{eq:def_L}
 \cL u(x)&=\int\limits_{\R^d\setminus \{0\}}\big(u(x+h)-u(x)-\langle
\nabla
u(x),h\rangle \mathbbm{1}_{B_1}(h)\big)K(x,h)\,dh \\
&=\frac12 \int\limits_{\R^d\setminus \{0\}}\big(u(x+h)-2u(x)+ u(x-h)
\big)K(x,h)\,dh \,,
\end{align}
under the assumption $K(x,h)=K(x,-h)$ and $K(x,h) \asymp |h|^{-d-\alpha}$ for
all $x$ and $h$ where
$\alpha \in (0,2)$ is fixed. Note that this class includes the case $\cL u =
-(-\Delta)^{\alpha/2} u$ and versions with bounded measurable coefficients.
As \cite{KrSa79} does, the article \cite{BL1} establishes a priori estimates in
Hölder spaces. Results like \autoref{lem:growth} have been obtained for
operators in the case $K(x,h) \asymp
|h|^{-d-\alpha}$ also for nonlinear problems, cf. \cite{Sil06}, \cite{CaSi09}
and \cite{GuSc12}.

The starting point of our research is the observation that \autoref{prop:hitting-prop} fails to hold for several interesting cases. One example
is given by $\cL$ as in \eqref{eq:def_L} with 
$K(x,h)=k(h) \asymp |h|^{-d}$ for $|h|\leq 1$ and some appropriate condition for
$|h| > 1$. For example, the geometric stable process with its generator
$-\ln(1+(-\Delta)^{\alpha/2})$, $0 < \alpha \leq 2$, can be represented
by \eqref{eq:def_L} with a kernel $K(x,h)=k(h)$ with such a behaviour for $|h|$
close to zero. The operator resp. the corresponding stochastic process can be
shown not to satisfy a uniformly
hitting estimate like \eqref{eq:hitting-estim}. This leads to the question
whether a priori estimates can be obtained by this approach at all.

Given a linear operator with bounded measurable coefficients of the form
\eqref{eq:def_L}, the main idea of this article is to determine an intrinsic
scale which allows to establish a modification of \eqref{eq:hitting-estim}.
We choose a measure different from the Lebesgue measure for the assumption
$|(B_{2R} \!\setminus\! B_R) \cap A| \geq \tfrac12 |B_{2R} \!\setminus\! B_R|$.

Let us formulate our assumptions and results. Assume $0 \leq \alpha < 2$ and
let $K\colon \R^d\times(\R^d\setminus\{0\})\rightarrow [0,\infty)$ be a
measurable function satisfying the following conditions:
{\allowdisplaybreaks
\Links
\begin{align}
\label{eq:K1}  \tag{$K_1$} & \ \ \ \sup\limits_{x\in
\R^d}\int\limits_{\R^d\setminus \{0\}}(1\wedge|h|^2)K(x,h)\,dh \leq
K_0 \,,\\
\label{eq:K2}  \tag{$K_2$} & \ \ \  K(x,h)=K(x,-h) \qquad (x\in
\R^d,\, h\in \R^d) \,,\\
\nonumber \\
 \label{eq:K3} \tag{$K_3$} & \ \ \  \kappa^{-1} \, \frac{\ell(|h|)}{|h|^d}\leq
K(x,h)\leq
\kappa \, \frac{\ell(|h|)}{|h|^d} \qquad (0<|h|\leq
1)
\end{align}
}for some numbers $K_0 > 0$, $\kappa>1$ and some function $\ell\colon
(0,1)\rightarrow (0,\infty)$ which is locally bounded and varies
regularly at zero with index $-\alpha\in (-2,0]$. Possible examples could be
$\ell(s)=1$,
$\ell(s)=s^{-3/2}$ and $\ell(s)=s^{-\beta} \ln(\tfrac{2}{s})^{2}$ for some
$\beta \in (0,2)$, see \autoref{sec:appendix} for a more detailed discussion.  

\pagebreak 

Suppose that there exists a strong Markov process $X=(X_t,\P_x)$ with
trajectories that are right continous with left limits associated with
${\cL}$ in the sense that for every $x \in \R^d$
\begin{itemize}
\item[(i)] $\P_x(X_0=x)=1$;
\item[(ii)] for any $f\in C_b^2(\R^d)$  the process $
\big\{f(X_t)-f(X_0)-\int_0^t
{\cL}f(X_s)\, ds |\, t \geq 0\big\}$ is a martingale under $\P_x$.
\end{itemize}

Note that the existence of such a Markov process comes for free in the
case when $K(x,h)$ is independent of $x$, see \autoref{sec:levystuff}. In
the general case it has been established by many authors in different general
contexts, see the discussion in \cite{AbKa09}. Denote by $\tau_A=\inf\{t>0|\,
X_t\not\in A\}$,
$T_A=\inf\{t>0|\, X_t \in A\}$  the first exit time resp. hitting time
of the
process $X$ for a measurable set $A\subset \R^d$.

\begin{Def}
 A bounded function $u\colon \R^d\rightarrow \R$ is said to be harmonic in an
open subset $D\subset \R^d$ with respect to $X$ (and ${\cL}$) if for
any bounded open set  $B\subset \overline{B}\subset D$ the stochastic process
$\{u(X_{\tau_B\wedge t})|\, t\geq 0\}$ is a $\P_x$-martingale for every $x\in
\R^d$\,.
\end{Def}

Before we can formulate our results we need to introduce an additional quantity.
Note
that \eqref{eq:K1} and \eqref{eq:K3} imply that $\int_0^1 s \,
\ell(s) \, \d s \leq c$ holds for some constant $c>0$. Let
$L\colon (0,1)\rightarrow (0,\infty)$ be defined by
$ L(r)=\int_r^1\frac{\ell(s)}{s}\, \d s$. The function $L$ is well
defined because $L(r) \leq r^{-2} \int_r^1 s^2 \frac{\ell(s)}{s}\, \d s \leq c
r^{-2}$.  See \autoref{sec:appendix} for several examples.  We note that
the function $L$ is always decreasing. Our main result concerning regularity is
the following result:

\begin{Thm}\label{thm:main} There exist constants $c>0$ and $\gamma\in (0,1)$
so that for all $r\in (0,\frac{1}{2})$ and $x_0\in \R^d$
\begin{align}\label{eq:main_reg}
|u(x)-u(y)|\leq c\|u\|_\infty
\frac{L(|x-y|)^{-\gamma}}{L(r)^{-\gamma}},\ \ x,y\in
B_{r/4}(x_0)
\end{align}
for all bounded functions $u\colon \R^d\rightarrow \R$ that are harmonic in
$B_r(x_0)$ with respect to $\cL$.
\end{Thm}

Let us comment on this result. It is important to note that the result
trivially holds if the function $L:(0,1) \to (0,\infty)$ satisfies
$\lim\limits_{r\to 0+}L(r) < +\infty$. This is equivalent to the condition
\begin{align} \label{eq:levy-finite}
\int\limits_{B_1} \frac{\ell(|h|)}{|h|^d} \, \d h < +\infty \,,
\end{align}
which, in the case $K(x,h) = k(h)$, means that the L\'{e}vy measure is finite. 
Thus, for the proof, we can concentrate on cases where \eqref{eq:levy-finite}
does not hold. One feature of this article is that our result holds true up to
and across the phase boundary determined by whether the kernel
$K(x,\cdot)$ is
integrable (finite L\'{e}vy measure) or not. 

Furthermore, note that the main result of \cite{BL1} is implied by
\autoref{thm:main} since
the choice $\ell(s)=s^{-\alpha}$, $\alpha \in (0,2)$, leads to
$L(r)\asymp r^{-\alpha}$. Given the whole spectrum of possible operators covered
by
our approach, this choice is a very specific one. It allows to use scaling
methods in the usual way which are not at our disposal here.
\autoref{tab:choices_l} in \autoref{sec:appendix} contains several admissible
examples one of which leads to
$L(0) < +\infty$ which means, as explained above,  that \eqref{eq:main_reg}
becomes pointless.

The main ingredient in the proof of \autoref{thm:main} is a new version of
\autoref{prop:hitting-prop} which we provide now. For $r\in (0,1)$ we
define a measure $\mu_r$ by 
\begin{equation}\label{eq:measure_mu_null}
 \mu_{r}(dx)=\frac{\ell(|x|)}{L(|x|) |x|^{d}} \, \mathbbm{1}_{B_1
\!\setminus\! B_r}(x) \,dx\,.
\end{equation}
Moreover, for $a>1$, we define a function $\varphi_a:(0,1) \to (0,1)$ by
$\varphi_a(r)=L^{-1}(\frac{1}{a}L(r))$. The following result is our
modification of \autoref{prop:hitting-prop}. 

\begin{Prop}\label{prop:hitting_new}
There exists a constant $c>0$ such that for all $a>1$, $r\in
(0,\frac{1}{2})$ and measurable sets $A\subset B_{\varphi_a(r)}\!\setminus\!
B_r$ with $\mu_{r}(A)\geq
\frac{1}{2}\mu_{r}(B_{\varphi_a(r)}\!\setminus\! B_r)$
 \[
  \P_x(T_A<\tau_{B_{\varphi_a(r)}}) \geq \P_x(X_{\tau_{B_r}}\in A)\geq
c \, \tfrac{\ln{a}}{a}
 \]
holds true for all $x\in B_{r/2}$.
\end{Prop}

The main novelties of \autoref{prop:hitting_new} are that the measure
$\mu_r$ depends on $r$ and that its density carries the factor $|x|^{-d}$.
These two changes allow us to deal with the classical cases as well as
with critical cases, e.g. given by $K(x,h)\asymp |h|^{-d}
\mathbbm{1}_{B_1}(h)$. 

The article is organised as follows: In \autoref{sec:levystuff} we review
the relation between translation invariant nonlocal operators and
semigroups/L\'{e}vy processes. Presumably, \autoref{prop:relation} is
interesting to many readers since it establishes a
one-to-one relation between the behavior of a L\'{e}vy measure at zero and the
multiplier of the corresponding generator for large values of $|\xi|$. In
\autoref{sec:prob_estim} we establish all tools needed to prove \autoref{prop:hitting_new} which is a special case of 
\autoref{prop:hitting_general}. \autoref{sec:regularity} contains the proof of
\autoref{thm:main}. The last section is \autoref{sec:appendix} in
which we collect important properties of regularly resp. slowly varying
functions. Moreover, the appendix contains a table with six examples which
illustrate the range of applicability of our approach.

Throughout the paper we use the notation $f(r)\asymp g(r)$ to denote
that the ration $f(r)/g(r)$ stays between two positive constants as $r$
converges to some value of interest.

\section{Translation invariant operators}\label{sec:levystuff}

The aim of this section is to discuss properties of the operator $\cL$
from \eqref{eq:def_L} in the translation invariant case, i.e. when 
$K(x,h)$ does not depend on $x \in \R^d$. In this case there is a one-to-one
correspondence between $\cL$ and multipliers, semigroups and stochastic
processes. One aim is to prove how the behavior of $\ell(|h|)$ for small values
of $|h|$ translates into properties of the multiplier or characteristic
exponent $\psi(|\xi|)$ for large values of $|\xi|$. This is acheived in
\autoref{prop:relation}. We add a subsection where we discuss which
regularity results are known in critical cases of the (much simpler) translation
invariant case. Note that our set-up, although
allowing for a irregular dependence of $K(x,h)$ on $x \in \R^d$, leads to new
results in these critical cases. 

\subsection{Generators of convolution semigroups and L\'{e}vy processes}

In this section we consider space homogeneous kernels of the form $K(x,h)=k(h)$
satisfing \eqref{eq:K1}--\eqref{eq:K3}. 
As we will see, the underlying stochastic process
belongs to the class of  L\' evy
processes\,. 

A stochastic process $X=(X_t)_{t\geq 0}$ on a probability space
$(\Omega,\mathcal{F},\P)$ is called a L\' evy process if it has
stationary and independent increments, $\P(X_0=0)=1$  and its
paths are $\P$-a.s. right continous with left limits\,. For $x\in \R^d$ we
define a $\P_x$
to be the law of the process $X+x$\,. In particular, $\P_x(X_t\in B)=\P(X_t\in
B-x)$ for $t\geq 0$ and measurable sets $B\subset \R^d$\,.

Due to stationarity and independence of increments, the characteristic function
of $X_t$ is given by
\[
 \E[e^{i\langle\xi, X_t\rangle}]=e^{-t\psi(\xi)},
\]
where $\psi$ is called characteristic exponent of $X$. It has the
following L\'
evy-Khintchine representation
\begin{equation}\label{eq:lk}
 \psi(\xi)=\frac{1}{2}\langle A\xi,\xi\rangle+\langle b,
\xi\rangle+\int_{\R^d\setminus
\{0\}}(1-e^{i\langle \xi,
h\rangle}+i\langle \xi,h\rangle \mathbbm{1}_{B_1}(h))\nu(dh)\,,
\end{equation}
where $A$ is a symmetric non-negative definite matrix , $b\in \R^d$ and $\nu$ is
a measure on
$\R^d\setminus \{0\}$ satisfying $\int_{\R^d\setminus\{0\}}(1\wedge
|y|^2)\nu(dy)<\infty$ called the
L\' evy measure of $X$. 

The converse also holds; that is, given $\psi$ as in the  L\' evy-Khintchine
representation \eqref{eq:lk}, there exists a L\' evy process $X=\{X_t\}_{t\geq
0}$ with the characteristic exponent $\psi$\,.
Details about L\' evy processes can be found in 
\cite{Be,S}\,.

To make a connection with our set-up, let $\nu$ be a measure defined by
$\nu(dh)=k(h)\,dh$. It follows from \eqref{eq:K1}--\eqref{eq:K3}  that $\nu$ is
a symmetric L\' evy measure. Let $X=\{X_t\}_{t\geq 0}$ be a L\' evy process
corresponding to the characteristic exponent $\psi$ as in \eqref{eq:lk} with
$A=0$, $b=0$ and the L\' evy measure $\nu(dh)=k(h)\,dh$\,.

Now, $P_t f(x):=\E_x[f(X_t)]$ defines a strongly continuous contraction semigoup
of operators $(P_t)_{t\geq 0}$ on the space $L^\infty(\R^d)$ equipped with
the essential-supremum norm.  Moreover, it is a
convolution semigroup, since 
\[
 \P_t f(x)=\E_0[f(x+X_t)]=\int_{\R^d}f(x+y)\mu_t(dy)\,,
\]
where $(\mu_t)_{t\geq 0}$ is a convolution semigroup of (probability) measures
defined by $\mu_t(B):=\P(X_t\in B)$.

The infinitesimal generator ${\cL}$ of the semigroup $(P_t)_{t\geq 0}$ is
given by
\begin{align}\label{eq:generator_levy}
 {\cL}u(x)=\int_{\R^d\setminus\{0\}}\big(u(x+h)-u(x)-\langle \nabla
u(x),h\rangle \mathbbm{1}_{B_1}(h)\big)k(h)\,dh
\end{align}
(cf. proof of \cite[Theorem 31.5]{S}). 

Since $\left\{u(X_t)-u(X_0)-\int_0^t{\cL}u(X_s)\,ds: t\geq 0\right\}$ is a
martingale (with respect to the natural filtration) for every $u\in C_b^2(\R^d)$
(cf. proof of \cite[Proposition VII.1.6]{RY}), it follows that $X$ is the
process which corresponds to the kernel $K(x,h)=k(h)$ in our set-up. 

It is worth of mentioning that there is a connection between the  characteristic
exponent and the symbol of the operator ${\cL}$. To be more precise, if
$\hat{f}(\xi)=\int_{\R^d}e^{i\xi\cdot x} f(x)\,dx$ denotes the Fourier transform
of a function $f\in L^1(\R^d)$, then 
\[
 \widehat{{\cL}f}(\xi)=-\psi(-\xi)\hat{f}(\xi)
\]
for any $f\in \mathcal{S}(\R^d)$, where $\mathcal{S}(\R^d)$ is the Schwartz
space (cf. \cite[Proposition I.2.9]{Be}). Hence $-\psi(-\xi)$ is the
symbol (multiplier) of the operator ${\cL}$\,.

We finish this section with the result that reveals connection between the
characteristic exponent $\psi$ and the function $L$\,.
\begin{Prop}\label{prop:relation} Let $\cL:\cS \to \cS$ be given by
\eqref{eq:generator_levy}. Assume $K(x,h):=k(h)$ satisfies
\eqref{eq:K1}-\eqref{eq:K3}. There is a constant $c>0$ such that 
\[
 c^{-1}L(|\xi|^{-1})\leq \psi(\xi)\leq cL(|\xi|^{-1}) \quad \text{ for } \xi\in
\R^d,\
|\xi|\geq 5\,.
\]
\end{Prop}

\proof
Note first that, by \eqref{eq:K3},
\begin{equation*}
 \kappa^{-1}j(|h|)\leq k(h)\leq \kappa j(|h|),\quad |h|\leq 1\,, 
\end{equation*}
where $j(s):=s^{-d}\ell(s)\,,\  s\in (0,1)$\,. 

Since $1-\cos{x}\leq \frac{1}{2}x^2$, it follows from
\eqref{eq:K1} and \eqref{eq:K3} that 
\begin{align*}
 \psi(\xi)&\leq \tfrac{1}{2}|\xi|^{2}\int_{|h|\leq
|\xi|^{-1}}|h|^2j(|h|)\,dh+2\int_{ |\xi|^{-1}<
|h|\leq 1}j(|h|)\,dh+2\int_{|h|>1}j(|h|)\,dh\\
&\leq c_1\left[|\xi|^2\int_0^{|\xi|^{-1}}s\ell(s)\,ds+L(|\xi|^{-1})+1\right]\\
&\leq c_2 (\ell(|\xi|^{-1})+L(|\xi|^{-1}))\leq c_3 L(|\xi|^{-1})\,,
\end{align*}
where in the first integral of the penultimate inequality Karamata's
theorem has
been used, while in the last inequality we have used that $\ell(s)\leq c_3L(s)$
for $s\in (0,1)$, cf. property \eqref{rem:l_versus_L} in
\autoref{sec:appendix}.

To prove the lower bound first we choose an orthogonal transformation of the
form $Oe_1=|\xi|^{-1}\xi$,
where $e_1:=(1,0,\ldots,0)\in \R^d$. Then a change of variable 
yields 
\begin{align*}
 \psi(\xi)&=\int_{\R^d\setminus\{0\}}(1-\cos(\xi\cdot
h))j(|h|)\,dh=\int_{\R^d\setminus\{0\}}(1-\cos{(|\xi|h_1)})j(|h|)\,dh\\
&\geq
\int_{[-1,1]^d}(1-\cos{(|\xi|h_1)})j(|h|)\,dh
\end{align*}
By the Fubini theorem,
\[
 \psi(\xi)\geq 2\int_0^1(1-\cos{(|\xi|r)})F(r)\,dr,
\]
where $
F(r):=\int_{[-1,1]^{d-1}}j(\sqrt{|z|^2+r^2})\,dz,\quad r\in
(0,\tfrac{1}{2})$\,. It follows from Potter's theorem (cf. property
\eqref{rem:potter} in \autoref{sec:appendix}) that there
is
a constant $c_4>0$ so that $j(r)\geq c_4 j(s)$ for all $0<r\leq s<1$. This
implies
\[
 F(r)\geq c_4 F(s),\qquad 0<r\leq s<1\,.
\]

Hence,
\begin{align*}
 \psi(\xi)&\geq 2 \sum_{k=0}^{\lfloor
\frac{\pi^{-1}|\xi|-\frac{3}{2}}{2}\rfloor}
\int_{|\xi|^{-1}(\frac{\pi}{2}+2k\pi)}^{|\xi|^{-1}(\frac{3\pi}{2}+2k\pi)}(1-\cos
{(|\xi|r)})F(r)\,
dr
\geq \frac{c_4\pi}{|\xi|}\sum_{k=0}^{\lfloor
\frac{\pi^{-1}|\xi|-\frac{3}{2}}{2}\rfloor}
F(|\xi|^{-1}(\tfrac{3\pi}{2}+2k\pi))\\
&\geq c_4^2\sum_{k=0}^{\lfloor \frac{\pi^{-1}|\xi|-\frac{3}{2}}{2}\rfloor}
\int_{|\xi|^{-1}(\frac{3\pi}{2}+2k\pi)}^{|\xi|^{-1}(\frac{3\pi}{2}+(2k+1)\pi)}
F(r)\,dr
\geq c_4^2\int_{\frac{3\pi}{2}|\xi|^{-1}}^1 F(r)\,dr\\
&\geq c_5 \int_{\frac{3\pi}{2}|\xi|^{-1}\leq |h|\leq 1} j(|h|)\,dh=c_6
L(\tfrac{3\pi}{2}|\xi|^{-1})\geq c_7 L(|\xi|^{-1})\,,
\end{align*}
where, in the last inequality, we have used property \eqref{rem:potter} from
\autoref{sec:appendix}. Note that
\cite{GR3} uses a similar trick to bound $\psi$ from below. 
\qed

\subsection{Known results in the translation invariant case}
Let us explain which results, related to Theorem \ref{thm:main}, have been
obtained in the case where $K(x,h)$ is independent of $x \in \R^d$.

H\"{o}lder estimates of harmonic functions are obtained for the L\' evy process
with the characteristic exponent
$\psi(\xi)=\frac{|\xi|^2}{\ln(1+|\xi|^2)}-1$ in \cite{Mi2} by establishing
a Krylov-Safonov type estimate replacing the Lebesgue measure with the capacity
of the sets involved. Recently, regularity estimates have been
obtained in \cite{GR3} for a class of isotropic unimodal L\' evy processes which
is quite general but does not include L\' evy processes with slowly varying L\'
evy exponents such as geometric stable processes. Regularity of harmonic
functions for such processes is investigated in \cite{Mi3}, where it is shown
that a result like \autoref{prop:hitting-prop} fails. Using the Green
function, logarithmic bounds for the modulus of continuity are obtained. At this
point it is worth mentioning that the transition density $p_t(x,y)$ of the 
geometric stable process satisfies $p_1(x,x)=\infty$, cf. \cite{SSV06}. This
illustrates that regularity results like \autoref{thm:main} in the case
$\ell(s)=1$ are quite delicate.

\section{Probabilistic estimates}\label{sec:prob_estim}
\begin{Prop}\label{prop:estp} There exists a constant $C_1>0$ such that
for $x_0\in \R^d$, $r\in
(0,1)$ and $t>0$ 
\[
\P_{x_0}(\tau_{B_r(x_0)}\leq t)\leq C_1t\,L(r) \,.
\]
\end{Prop}
\proof
Let $x_0\in \R^d$, $0<r<1$ and let $f\in C^2(\R^d)$ be a positive function such that
\[
f(x)=\left\{\begin{array}{cl}
|x-x_0|^2, & |x-x_0|\leq \frac{r}{2}\\
r^2,& |x-x_0|\geq r
\end{array}\right.
\]
and for some $c_1 >0$ 
\[
|f(x)|\leq c_1 r^2, \ \ \left|\frac{\partial f}{\partial x_i}(x)\right|\leq c_1 r \ \ \textrm{ and }
\ \ \left|\frac{\partial^2 f}{\partial x_i \partial x_j}(x)\right|\leq c_1.
\]
By the optional stopping theorem we get
\begin{align}\label{eq:optional_stopping}
\E_x f(X_{t\wedge \tau_{B_r(x_0)}})-f(x_0)=\E^x \int_0^{t\wedge
\tau_{B_r(x_0)}}\mathcal{L}f(X_s)\,
ds,\ \ t>0.
\end{align}
Let $x\in B_r(x_0)$. We estimate $\mathcal{L}f(x)$ by splitting the integral in
\eqref{eq:def_L} into three parts.
\begin{align*}
\int_{B_r} &(f(x+h)-f(x)-\nabla f(x)\cdot h \mathbbm{1}_{\{|h|\leq
1\}})K(x,h)\,dh\\
& \leq c_2\int_{B_r}|h|^2K(x,h)\,dh\leq c_2\kappa
\int_{B_r}|h|^{2-d}\ell(|h|)\,dh \leq c_3 r^{2}\ell(r),
\end{align*}
where in the last line we have used Karamata's theorem, cf. property
\eqref{rem:karamata} in \autoref{sec:appendix}. On the other hand, on $B_r^c$ we
have
\begin{align*}
\int_{B_r^c} &(f(x+h)-f(x))K(x,h)\,dh\leq 2\|f\|_\infty
\int_{B_r^c}K(x,h)\,dh\\
& \leq 2\|f\|_\infty \left(\kappa\int_{B_1\setminus
B_r}|h|^{-d}\ell(|h|)\,dh+\int_{B_1^c}K(x,h)\,dh\right) \leq c_4r^{2}L(r)\,dr
\,,
\end{align*}
where we applied property \eqref{rem:L_notzero} from \autoref{sec:appendix}.
Last, we estimate
\begin{align*}
\left|\int_{B_1\setminus B_r}h\cdot \nabla f(x) K(x,h)\,dh\right|&\leq
c_1 r\int_{B_1\setminus B_r} |h| K(x,h)\,dh\\
&\leq c_1 \kappa r \int _
{B_1\setminus B_r} |h|^{-d+1}\ell(|h|)\,dh \leq c_5 r^2\ell(r),
\end{align*}
by Karamata's theorem again. Therefore, by property \eqref{rem:l_versus_L} from
\autoref{sec:appendix} we conclude that
there is a constant $c_6>0$ such that for all $x\in B_r(x_0)$ and $r\in (0,1)$
we have
\begin{equation}\label{eq:prop1e1}
\mathcal{L}f(x)\leq c_6 r^2 L(r).
\end{equation}
Let us look again at \eqref{eq:optional_stopping}. On
$\{\tau_{B_r(x_0)}\leq t\}$ we have $X_{t\wedge \tau_{B_r(x_0)}}\in
B_r(x_0)^c$ and so
$f(X_{t\wedge \tau_{B_r(x_0)}})\geq r^2$. Thus, by  (\ref{eq:prop1e1}) and
\eqref{eq:optional_stopping} we get
\[
\P_{x_0}(\tau_{B_r(x_0)}\leq t)\leq c_6 L(r)t.
\]
\qed

\begin{Prop}\label{prop:prob-2}
	There are constants $C_2>0$ and $C_3>0$ such that for $x_0\in \R^d$
	\[
		\sup_{x\in \R^d}\E_x\tau_{B_r(x_0)}\leq
\frac{C_2}{L(r)}\,,\quad r\in
(0,1/2)
	\]
and
	\[
		\inf_{x\in B_{r/2}(x_0)}\E_x\tau_{B_r(x_0)}\geq
\frac{C_3}{L(r)}\,,\quad r\in
(0,1)
	\]

\end{Prop}
\proof
	The proof is similar to the proof of the exit time estimates in \cite{BL1}. 

      (a) First we prove the upper estimate for the exit time.
	Let $x\in\R^d$, $r\in (0,1/2)$ and let
	\[
		S=\inf\{t>0|\, |X_t-X_{t-}|>2r\}
	\]
	be the first time of a jump larger than $2r$. With the help of the
L\' evy system formula (cf. \cite[Proposition 2.3]{BL1}) and \eqref{eq:K3} we
can deduce
{\allowdisplaybreaks 
	\begin{align}
		\P_x(S\leq L(r)^{-1}) & =\E_x\sum_{t\leq L(r)^{-1}\wedge S}
\mathbbm{1}_{\{|X_t-X_{t-}|>2r\}} =\E_x\int\limits_0^{L(r)^{-1}\wedge
S}\int\limits_{B_{2r}^c}K(X_s,h)\,dh\,ds \nonumber\\
		&\geq c_1\E_x[L(r)^{-1}\wedge S]\int\limits_{2r}^1\frac{\ell(t)}{t}\,dt\,.\label{eq:prob-6}
	\end{align}
	}
	Since $L$ is regularly varying at zero, 
	\begin{align*}
		\E_x[L(r)^{-1}\wedge S]& \geq L(r)^{-1}\P_x(S>L(r)^{-1}) \geq
c_2L(2r)^{-1}\big(1-\P_x(S\leq L(r)^{-1})\big)\,
	\end{align*}
	and so it follows from (\ref{eq:prob-6}) that
	\begin{equation}\label{eq:prob-7}
		\P_x(S\leq L(r)^{-1})\geq c_3
	\end{equation}
	with $c_3=\frac{c_1c_2}{c_1c_2+1} \in (0,1)$. The strong Markov
property and (\ref{eq:prob-6}) lead to
	\[
		\P_x(S>m L(r)^{-1})\leq (1-c_3)^m,\ \ m\in \N\,.
	\]
	Since $\tau_{B_r(x_0)}\leq S$, 
	\begin{align*}
		\E_x\tau_{B_r(x_0)}\leq \E_x S &\leq
L(r)^{-1}\sum_{m=0}^\infty (m+1)\P_x(S>L(r)^{-1}m)\\
		&\leq L(r)^{-1}\sum_{m=0}^\infty (m+1)(1-c_3)^m\,.
	\end{align*}

(b) Now we prove the lower estimate of the exit time. Let $r\in
(0,1)$ and
$y\in B_{r/2}(x_0)$.
By \autoref{prop:estp},
\[
  \P_y(\tau_{B_r(x_0)}\leq t)\leq \P_y(\tau_{B_{r/2}(y)}\leq t)\leq
C_1tL(r/2),\quad t>0\,,
\]
since $B_{r/2}(y)\subset B_r(x_0)$\,. Choose $t=\frac{1}{2C_1L(r/2)}$. Then
\begin{align*}
 \E_y\tau_{B_r(x_0)}&\geq \E_y[\tau_{B_r(x_0)};\tau_{B_r(x_0)}>t]\geq
t\P_y(\tau_{B_r(x_0)}>t)\\
		    &\geq t(1-C_1L(r/2)t)=\frac{1}{4C_1L(r/2)}\,.
\end{align*}
By \eqref{rem:l-reg_L-reg} from \autoref{sec:appendix} we know that
$L$ is regularly varying at zero. Hence there is a
constant $c_1>0$ such that $L(r/2)\leq c_1 L(r)$ for all $r\in (0,1/2)$.
Therefore
$\E_y\tau_{B_r(x_0)}\geq \frac{1}{4C_1c_1L(r)}$\,.

\qed

\begin{Prop}\label{prop:prob-3}
	There is a constant $C_4>0$ such that for all $x_0\in \R^d$ and $r,s\in
(0,1)$ satisfying $2r<s$
	\[
		\sup_{x\in B_r(x_0)}\P_x(X_{\tau_{B_r(x_0)}}\not\in
B_s(x_0))\leq C_4 \frac{L(s)}{L(r)}\,.
	\]
\end{Prop}
\proof
	Let $x_0 \in \R^d$, $r,s \in (0,1)$ and $x\in B_r(x_0)$. Set
$B_r:=B_r(x_0)$. By the L\' evy system formula, for $t>0$
	\begin{align*}
		\P_x(X_{\tau_{B_r}\wedge t}\not\in B_s)&=\E_x\sum\limits_{s\leq
\tau_{B_r}\wedge t} \mathbbm{1}_{\{X_{s-}\in B_r,X_s\in
B_s^c\}} =\E_x\int\limits_0^{\tau_{B_r}\wedge
t}\int\limits_{B_s^c} K(X_s,z-X_s)\,dz\,ds \,.
	\end{align*}
	
	Let  $y\in B_r$. Since $s\geq 2r$, it follows that
$B_{s/2}(y)\subset B_s$ and hence
	\begin{align*}
		\int\limits_{B_s^c} K(y,z-y)\,dz&\leq
\int\limits_{B_{s/2}(y)^c} K(y,z-y)\,dz \leq
c_1\int_{s/2}^1 \frac{\ell(u)}{u}\,du+c_2 \leq c_3L(s)\,.
	\end{align*}
	where in the last inequality we have used that $L$ varies
regularly at zero and that $\lim\limits_{r\to 0+}L(r)> 0$, cf.
\eqref{rem:L_notzero} in \autoref{sec:appendix}. 
	
	The above considerations together with \autoref{prop:prob-2} imply
	\[
		\P_x(X_{\tau_{B_r}\wedge t}\not\in B_s)\leq c_3 L(s)\E_x \tau_{B_r}\leq c_4 \frac{L(s)}{L(r)}\,.
	\]
	Letting $t\to\infty$ we obtain the desired estimate.
\qed

For $x_0\in \R^d$ and $r\in (0,1)$ we define the following measure
\begin{equation}\label{eq:measure_mu}
 \mu_{x_0,r}(dx)=\frac{\ell(|x-x_0|)}{L(|x-x_0|)}\,|x-x_0|^{-d}
\mathbbm{1}_{\{r\leq
|x-x_0|<1\}}\,dx\,.
\end{equation}

Define $\varphi_a(r)=L^{-1}(\frac{1}{a}L(r))$ for $r\in (0,1)$ and $a>1$. The
following property is important for the construction below:
\begin{align} 
r = L^{-1}(L(r)) \leq L^{-1}(\tfrac{1}{a}L(r)) = \varphi_a(r)
\,. 
\end{align}

Now we can prove a Krylov-Safonov type hitting estimate which
includes \autoref{prop:hitting_new} as a special case. 

\begin{Prop}\label{prop:hitting_general}
There exists a constant $C_5>0$ such that for all $x_0\in \R^d$, $a>1$, $r\in
(0,\frac{1}{2})$ and
$A\subset B_{\varphi_a(r)}(x_0)\setminus B_r(x_0)$ satisfying
$\mu_{x_0,r}(A)\geq
\frac{1}{2}\mu_{x_0,r}(B_{\varphi_a(r)}(x_0) \setminus
B_r(x_0))$
 \[
  \P_y(T_A<\tau_{B_{\varphi_a(r)}(x_0)})\geq \P_y(X_{\tau_{B_r(x_0)}}\in A)\geq
C_5\frac{\ln{a}}{a}\,,\quad y\in B_{r/2}(x_0)\,.
 \]
\end{Prop}

\proof
Consider $x_0\in \R^d$, $a>1$, $r\in
(0,\frac{1}{2})$ and a set 
$A\subset B_{\varphi_a(r)}(x_0)\setminus B_r(x_0)$ satisfying
$\mu_{x_0,r}(A)\geq
\frac{1}{2}\mu_{x_0,r}(B_{\varphi_a(r)}(x_0) \setminus
B_r(x_0))$. Set $\mu:=\mu_{x_0,r}$, $\varphi:=\varphi_a$,  $B_s:=B_s(x_0)$ and
let $y\in
B_{r/2}$. The first
inequality follows from 
$\{X_{\tau_{B_r}}\in A\}\subset \{T_A<\tau_{B_{\varphi(r)}}\}$ since $A\subset
B_{\varphi(r)}\setminus
B_r$\,.

By the L\' evy system formula, for $t>0$, 
	\begin{align}\label{eq:ks-0}
		\P_y(X_{\tau_{B_r}\wedge t}\in A)&=\E_y\sum\limits_{s\leq \tau_{B_r}\wedge
t} \mathbbm{1}_{\{X_{s-}\in B_r,X_s\in A\}}	
=\E_y\int\limits_0^{\tau_{B_r}\wedge t}\int\limits_{A} K(X_s,z-X_s)\,dz\,ds\,.
	\end{align}

Since $|z-x|\leq |z-x_0|+|x_0-x|\leq |z-x_0|+r\leq 2|z-x_0|$ for $x\in B_r$ and $z\in B_r^c$,
\begin{equation}\label{eq:ks-1}
 \E_y\int\limits_0^{\tau_{B_r}\wedge t}\int\limits_{A} K(X_s,z-X_s)\,dz\,ds\geq
c_1\E_y[\tau_{B_r}\wedge t]\int_A \frac{\ell(|z-x_0|)}{|z-x_0|^d}\,dz\,,
\end{equation}
where we have used property \eqref{rem:potter} given in
\autoref{sec:appendix}.

Since $L$ is decreasing,
\begin{align}
 \int_A \frac{\ell(|z-x_0|)}{|z-x_0|^d}\,dz&=\int_A
L(|z-x_0|)\mu(dz) \geq
L(\varphi(r))\mu(A)\geq
\frac{L(r)}{2a}\mu(B_{\varphi(r)}\setminus B_r)\,.\label{eq:ks-2}
\end{align}
Noting that 
\[\mu(B_{\varphi(r)}\setminus B_r)=c_2\int_r^{\varphi(r)}
\frac{1}{L(s)}\frac{\ell(s)\,ds}{s}=-c_2\ln L(s)|_r^{\varphi(r)}=c_2\ln a\,,\]
we conclude from
(\ref{eq:ks-0})--(\ref{eq:ks-2}) that
\[
  \P_y(T_A<\tau_{B_{\varphi_a(r)}(x_0)})\geq c_3
L(r)\frac{\ln{a}}{a}\E_y[\tau_{B_r}\wedge t]\,.
\]
Letting $t\to\infty $ and using the lower bound in \autoref{prop:prob-2} we get
\begin{align*}
 \P_y(T_A<\tau_{B_{\varphi_a(r)}(x_0)}) \geq c_3
L(r)\,\frac{\ln{a}}{a}\,\E_y\tau_{B_r} \geq
c_3L(r)\,\frac{\ln{a}}{a}\,C_3L(r)^{-1}=c_3C_3\frac{\ln{a}}{a}\,.
\end{align*}

\qed

\section{Reglarity of harmonic functions}\label{sec:regularity}

\proof[Proof of \autoref{thm:main}] 

Let $x_0 \in \R^d$, $r\in (0,\frac12)$, $x\in B_{r/4}(x_0)$. Using
\eqref{rem:potter} from \autoref{sec:appendix} with $\delta=1$, we see
that there is a constant
$c_0\geq 1$ so that 
\begin{align}\label{eq:potterp}
 \frac{L(s)}{L(s')}\leq c_0\left(\frac{s}{s'}\right)^{-\alpha-1},\quad
0<s<s'<1\,.
\end{align}

Define for $n\in \N$
\[
	r_n:=L^{-1}( L(\tfrac{r}{2})a^{n-1})\quad \text{ and }\quad s_n:=3 \|u\|_\infty
b^{-(n-1)}
\]
for some constants $b\in (1,\frac{3}{2})$ and $a>c_0 2^{\alpha+1}$ that will be
chosen in the proof
independently of $n$, $r$ and $u$. As we explained in the
introduction, \autoref{thm:main} trivially holds true of $\lim\limits_{r\to 0+}
L(r)$
is finite. Thus, we can assume $\lim\limits_{r\to 0+} L(r)$ to be infinite. This
implies that $r_n \to 0$ for $n \to \infty$ as it should be. 

We will use the following abbreviations:
\[
	B_n:=B_{r_n}(x), \quad \tau_n:=\tau_{B_n}, \quad m_n:=\inf_{B_n}u,
\quad M_n:=\sup_{B_n}u\,.
\]

We are going to prove 
 \begin{equation}\label{eq:thm-1}
 	M_k-m_k\leq s_k
 \end{equation}
 for all $k\geq 1$. 
 
Assume for a moment that \eqref{eq:thm-1} is proved. Then, for any $r\in
(0,\frac{1}{2})$ and $y\in B_{r/4}(x_0)\subset B_{r/2}(x)$ we can find $n\in
\N$ so that 
 \[
 	r_{n+1}\leq |y-x|< r_n\,.
 \]
Furthermore, since $L$ is decreasing, we obtain with 
$\gamma=\frac{\ln{b}}{\ln{a}}\in (0,1)$  
\begin{align*}
 	|u(y)-u(x)|&\leq s_n=3 b\|u\|_\infty a^{-n\frac{\ln b}{\ln a}} =3
b\|u\|_\infty \left[\frac{L(r_{n+1})}{
L(\frac{r}{2})}\right]^{-\frac{\ln{b}}{\ln{a}}} \leq 3
b\|u\|_\infty\left[\frac{L(|x-y|)}{L(\frac{r}{2})}\right]^{-\gamma}\,,
\end{align*}
which proves our assertion. Thus it remains to prove
(\ref{eq:thm-1}).

We are going to prove \eqref{eq:thm-1} by an inductive argument.
Obviously, $M_1-m_1\leq 2\|u\|_\infty\leq
s_1$. Since $1<b<\frac{3}{2}$, it follows that 
\[
M_2-m_2\leq 2\|u\|_\infty\leq 3\|u\|_\infty b^{-1}=s_2\,.
\]

Assume now that (\ref{eq:thm-1})  is true for all $k\in \{1,2,\ldots,n\}$ for some $n\geq 2$. 

Let $\varepsilon>0$ and take ${z_1},z_2 \in B_{n+1}$ so that
\[
	u({z_1})\leq m_{n+1}+\frac{\varepsilon}{2}\ \ \ \ \ \ \ \ u(z_2 )\geq
M_{n+1}-\frac{\varepsilon}{2}\,.
\]
It is enough to show that
\begin{equation}\label{eq:thm-2}
	u(z_2 )-u({z_1})\leq s_{n+1},
\end{equation}
since then 
\[
	M_{n+1}-m_{m+1}-\varepsilon\leq s_{n+1},
\]
which implies (\ref{eq:thm-1})  for $k=n+1$, since $\varepsilon>0$ was arbitrary. 

By the optional stopping theorem,
\begin{align*}
	u(z_2 )-u({z_1})=&\ \E_{z_2}[u(X_{\tau_n})-u({z_1})]\nonumber \\
	=&\ \E_{z_2}[u(X_{\tau_n})-u({z_1});X_{\tau_n}\in
B_{n-1}]\\&+\sum\limits_{i=1}^{n-2}\E_{z_2}[u(X_{\tau_n})-u({z_1});X_{\tau_n}\in
B_{n-i-1}\setminus B_{n-i}]\\&+\E_{z_2}[u(X_{\tau_n})-u({z_1});X_{\tau_n}\in
B_1^c]=I_1+I_2+I_3\,.
\end{align*}

Let $A=\{z \in B_{n-1}\setminus B_n|\, u(z)\leq \frac{m_n+M_n}{2}\}$. It is
sufficient to consider the case $\mu_{x,r_n}(A)\geq
\frac{1}{2}\mu_{x,r_n}(B_{n-1}\setminus B_n)$, where
$\mu_{x,r}$ is the measure defined by (\ref{eq:measure_mu}). In the
remaining case we would use 
$\mu_{x,r_n}((B_{n-1}\setminus B_n)\setminus A)\geq \frac{1}{2}\mu_{x,r_n}(B_{n-1}\setminus B_n)$
 and could continue the proof with  $\|u\|_\infty -u$ and 
 \[
 	(B_{n-1}\setminus B_n)\setminus A=\left\{z\in B_{n-1}\setminus B_n
|\, \|u\|_\infty
-u(z)\leq
\frac{\|u\|_\infty
-m_n+\|u\|_\infty -M_n}{2}\right\}
 \]
instead of $u$ and $A$.

The estimate \eqref{eq:potterp} implies  
$a=\tfrac{L(r_{n+1})}{L(r_n)}\leq c_0(\tfrac{r_{n+1}}{r_n})^{-\alpha-1}$, from
where we deduce $r_{n+1}\leq r_n (c_0a^{-1})^{\frac{1}{\alpha+1}} \leq
\frac{r_n}{2}$ because of $a>c_0 2^{\alpha+1}$. Next, we make use of the
following property:
\begin{align}\label{eq:role_varphi} 
r_{n-1} = L^{-1}(L(\tfrac{r}{2}) a^{n-2}) = L^{-1}(\tfrac{1}{a}
L(\tfrac{r}{2}) a^{n-1}) =  L^{-1}(\tfrac{1}{a} L(r_n)) = \varphi_a(r_n) \,. 
\end{align}
Then by
\autoref{prop:hitting_general} (with
$r=r_n$ and $x_0=x$) we get  \[p_n:=\P_{z_2}(X_{\tau_n}\in A)\geq
C_5\frac{\ln{a}}{a}\,.\]

Hence,
\begin{align*}
 I_1&=\E_{z_2}[u(X_{\tau_n})-u({z_1});X_{\tau_n}\in B_{n-1}]\\
&=\E_{z_2}[u(X_{\tau_n})-u({z_1});X_{\tau_n}\in
A]+\E_{z_2}[u(X_{\tau_n})-u({z_1});X_{\tau_n}\in B_{n-1}\setminus
A]\\
&\leq \left(\tfrac{m_n+M_n}{2}-m_n\right)p_n+s_{n-1}(1-p_n)\\
&\leq \tfrac{1}{2}s_n p_n+s_{n-1}(1-p_n)\leq s_{n-1}(1-\tfrac{1}{2}p_n)\leq
s_{n-1}(1-\tfrac{C_5\ln{a}}{2a})\,.
\end{align*}

By \autoref{prop:prob-3}, 
\begin{align*}
	I_2& \leq \sum\limits_{i=1}^{n-2} s_{n-i-1}\P_{z_2}(X_{\tau_n}\not\in
B_{n-i})
	\leq C_4 \sum\limits_{i=1}^{n-2} s_{n-i-1}\tfrac{L(r_{n-i})}{L(r_n)}\\
	&\leq 3C_4\|u\|_\infty  \sum\limits_{i=1}^{n-2}
b^{-(n-i-2)}\tfrac{a^{n-i-1}}{a^{n-1}}\leq
3C_4 \|u\|_\infty \tfrac{b^{-n+3}}{a-b}\\
	&\leq C_4  \tfrac{b^3}{a-b}s_{n+1}\,.
\end{align*}
Similarly, by \autoref{prop:prob-3}, 
\[
	I_3\leq 2\|u\|_\infty\P_{z_2}(X_{\tau_n}\not\in B_1)\leq 2C_4
\|u\|_\infty
\tfrac{L(r_1)}{L(r_n)}=\tfrac{2C_4 }{3}
b\left(\tfrac{b}{a}\right)^{n-1}s_{n+1}\leq C_4 
\tfrac{b^2}{a} s_{n+1}\,.
\]
Hence,
\[
u(z_2)-u(z_1)\leq s_{n+1}b^2\left[1-\tfrac{C_5\ln{a}}{2a}+ \tfrac{C_4 b}{a-b}+
\tfrac{C_4 }{a}\right]\,.
\]

Since $a-b\geq \frac{a}{4}$ for $b\in (1,\frac{3}{2})$ and
$a>c_0 2^{\alpha+1}\geq 2$, it follows
that
\[q:=1-\tfrac{C_5\ln{a}}{2a}+
\tfrac{C_4  b}{a-b}+
\tfrac{C_4 }{a}\leq 1-\tfrac{C_5\ln{a}}{2a}+\tfrac{7C_4 }{a}=1-\tfrac{C_5\ln
a-14C_4 }{2a}\,.\]

Next, we choose $a>c_0 2^{\alpha+1}$ so large that $C_5\ln{a}-14C_4 >0$. Thus
$q<1$.
Finally, we choose $b\in
(1,\frac{3}{2})$ sufficiently small so that $b^2q<1$\,. 

Hence, (\ref{eq:thm-2}) holds, which finishes the proof of the inductive step and the theorem\,.
\qed

\begin{appendix}

\section{Slow and Regular Variation} \label{sec:appendix}

In this section we collect some properties of slowly resp. regularly varying
functions that are used in our main arguments. Moreover we list several examples
which illustrate the range of application of our approach.

\begin{Def}
 A measurable and positive function $\ell\colon (0,1)\rightarrow (0,\infty)$ is
said to vary regularly at zero with index $\rho\in
\R$ if for every $\lambda>0$
\[
 \lim_{r\to 0+}\frac{\ell(\lambda r)}{\ell(r)}=\lambda^{\rho}\,.
\]
If a function varies regularly at zero with index $0$ it is said to vary
slowly at zero. For simplicity, we call such functions \emph{regularly
varying} resp. \emph{slowly varying} functions.
\end{Def}

Note that slowly resp. regularly varying functions include functions which are
neither
increasing nor decreasing. By \cite[Theorem 1.4.1 (iii)]{BGT} it follows that
any function $\ell$ that
varies regularly with index $\rho\in \R$ is of the form
$\ell(r)=r^{\rho}\ell_0(r)$ for some
function $\ell_0$ that varies slowly. 

Assume $\int_0^1 s \, \ell(s) \, \d s \leq c$ for some $c >0$. Let $L\colon
(0,1)\rightarrow (0,\infty)$ be defined by
\[
	L(r)=\int\limits_r^1\frac{\ell(s)}{s}\, \d s\ .
\]
The function $L$ is well defined because $L(r) = r^{-2} \int_r^1
r^2 \frac{\ell(s)}{s}\, \d s \leq r^{-2} \int_r^1
s \ell(s) \, \d s \leq c r^{-2}$. Note that \eqref{eq:K1} and
\eqref{eq:K3} imply that $\int_0^1 s \, \ell(s) \, \d s \leq c$ does hold in our
setting. We note that the function $L$ is always decreasing. 

Let us list further properties which are making use of in our proofs. Note that
they are established \cite{BGT} for functions which are slowly resp. regularly
varying at the point $+\infty$. By a simple inversion we adopt the results
to functions which are slowly resp. regularly
varying at the point $0$.

\begin{enumerate}
 \item \label{rem:l_versus_L} If $\ell$ is slowly varying, then
\cite[Proposition
1.5.9a]{BGT} $L$ is slowly varying  with 
\[\lim\limits_{r\to
0+}L(r)=+\infty\qquad \text{ and }\qquad \lim\limits_{r\to
0+}\frac{\ell(r)}{L(r)}=0\,.\]
\item \label{rem:karamata} If $\ell$ is slowly varying and $\rho>-1$, then
Karamata's theorem \cite[Proposition 1.5.8]{BGT} ensures 
\[
  \lim_{r\to 0+}\frac{\int_0^r
s^{\rho}\ell(s)\,ds}{r^{\rho+1}\ell(r)}=(\rho+1)^{-1}\,.
\]
\item \label{rem:l-reg_L-reg} If $\ell$ is
regularly varying of order $-\alpha < 0$ (in our case $0 < \alpha < 2$),
then
\cite[Theorem 1.5.11]{BGT} 
\[
 \lim_{r\to 0+}\frac{L(r)}{\ell(r)}=\alpha^{-1}\,.
\] In particular, if $\ell$ is regularly varying of order $-\alpha < 0$,
then so is $L$. 
\item \label{rem:potter} Assume $\ell$ is regularly varying of order
$-\alpha \leq 0$ and stays bounded away from $0$ and $+\infty$ on every compact
subset of $(0,1)$. Then Potter's theorem \cite[Theorem 1.5.6 (ii)]{BGT} implies
that for every $\delta>0$ there is a
constant
$C=C(\delta)\geq 1$ such that for $r,s \in (0,1)$ 
\begin{align*}
 \frac{\ell(r)}{\ell(s)}\leq C
\max\left\{ \left(\frac{r}{s}\right)^{-\alpha-\delta},
\left(\frac{r}{s}\right)^{-\alpha+\delta} \right\} \,.
\end{align*}
\item \label{rem:L_notzero} Since $L$ is nonincreasing, we observe
$\lim\limits_{r\to 0+} L(r) \in (0,+\infty]$.
\end{enumerate}

\smallskip

\renewcommand{\arraystretch}{1.9}

\begin{table}[ht]
\caption{Different choices for the function $\ell$ when $\beta \in (0,2)$, $a >
1$.}
\begin{tabular}{ c  c  c  c  }
\hline
{\bf No.\! (i)} & $\ell_i(s)$ &
$L_i(s)$ & $\varphi_a(s)=L_i^{-1}(\frac{1}{a} L_i(s))$ \\
\toprule
$1$ & $s^{-\beta} \, \ln(\frac{2}{s})^{2}$ & $\asymp s^{-\beta} \, \ln
(\frac{2}{s})^{2}$ & $\asymp 
s$\\
\hline
$2$ & $s^{-\beta}$ & $\frac{1}{\beta} (s^{-\beta}-1)$ & $\asymp 
s$\\
\hline
$3$ & $\ln(\frac{2}{s})$ & $\asymp \ln^2(\frac{2}{s})$ &
$\asymp s^{1/\sqrt{a}}$\\
\hline
$4$ & $1$ &  $\ln(\frac{1}{s})$ & $ s^{1/a}$\\
\hline
$5$ & $\ln(\frac{2}{s})^{-1}$ & $\asymp \ln(\ln(\frac{2}{s}))$ & $\asymp \exp(-
(\ln(\frac{2}{s}))^{1/a})$\\
\hline
$6$ & $\ln(\frac{2}{s})^{-2}$ & $\ln(2)^{-1}-\ln(\tfrac{2}{s})^{-1}$ &
$\asymp \exp(-(\frac{a-1}{a\ln(2)}+\frac{1}{a\ln(2/s)})^{-1})$ \\
\bottomrule
\end{tabular}
\label{tab:choices_l}
\end{table}


Let us look at different choices for the function $\ell$, given in Table
\ref{tab:choices_l}. Here  $\beta \in (0,2)$, $a >
1$ are fixed. We list six examples of a function $s \mapsto \ell_i(s)$ together
with $s\mapsto L_i(s)$ and $s \mapsto \varphi_a(s)=L_i^{-1}(\frac{1}{a}
L_i(s))$. Recall that the function $\varphi_a$ appears in
\autoref{prop:hitting_new} and determines the scaling that we are using, see
also property \eqref{eq:role_varphi} and the definition of $r_n$ in the proof of
\autoref{thm:main}. Note that case No. 6 is significantly different from the
other
cases. Both, the integral $\int_{B_1}|h|^{-d} \ell_6(|h|) \, \d h$ and the
expression $\lim\limits_{s\to 0+} L_6 (s)$ are finite. Moreover, the limit
$\lim\limits_{s\to
0+} L_6^{-1}(\frac{1}{a} L_6(s))$ is not equal to zero. These differences
reflect the
fact that the corresponding operator in \eqref{eq:def_L} has an integrable
kernel. Recall that
\autoref{prop:relation} relates the behavior of the function $L$ close
to the origin to the behaviour of the multiplier of the operator (in the case
of constant coefficents) for large values of $|\xi|$. In the case No. 6 the
multiplier stays bounded.

\end{appendix}


{\bf Acknowledgements:} We thank T. Grzywny for a helpful comment on the limit
case $\alpha=2$. 

\bibliography{gen_kernels}

\providecommand{\bysame}{\leavevmode\hbox to3em{\hrulefill}\thinspace}
\providecommand{\MR}{\relax\ifhmode\unskip\space\fi MR }
\providecommand{\MRhref}[2]{%
  \href{http://www.ams.org/mathscinet-getitem?mr=#1}{#2}
}
\providecommand{\href}[2]{#2}
\begin{thebibliography}{Mim13b}

\bibitem[AK09]{AbKa09}
H.~Abels and M.~Kassmann, \emph{The {C}auchy problem and the martingale problem
  for integro-differential operators with non-smooth kernels}, Osaka J. Math.
  \textbf{46} (2009), no.~3, 661--683. \MR{2583323 (2011d:35505)}

\bibitem[Ber96]{Be}
J.~Bertoin, \emph{L\' evy processes}, Cambridge University Press, Cambridge,
  1996.

\bibitem[BGT87]{BGT}
N.~H. Bingham, C.~M. Goldie, and J.~L. Teugels, \emph{Regular variation},
  Cambridge University Press, Cambridge, 1987.

\bibitem[BL02]{BL1}
R.~F. Bass and D.~Levin, \emph{Harnack inequalities for jump processes},
  Potential Anal. \textbf{17} (2002), 375--388.

\bibitem[CC95]{CaCa95}
Luis~A. Caffarelli and Xavier Cabr{\'e}, \emph{Fully nonlinear elliptic
  equations}, American Mathematical Society Colloquium Publications, vol.~43,
  American Mathematical Society, Providence, RI, 1995. \MR{1351007 (96h:35046)}

\bibitem[CS09]{CaSi09}
L.~Caffarelli and L.~Silvestre, \emph{Regularity theory for fully nonlinear
  integro-differential equations}, Comm. Pure Appl. Math. \textbf{62} (2009),
  no.~5, 597--638. \MR{2494809 (2010d:35376)}

\bibitem[DGV12]{DBGV12}
E.~DiBenedetto, U.~Gianazza, and V.~Vespri, \emph{Harnack's inequality for
  degenerate and singular parabolic equations}, Springer Monographs in
  Mathematics, Springer, New York, 2012. \MR{2865434}

\bibitem[Grz13]{GR3}
T.~Grzywny, \emph{On {H}arnack inequality and {H}{\"{o}}lder regularity for
  isotropic unimodal {L}\'{e}vy processes}, Potential Anal. (2013), to appear.

\bibitem[GS12]{GuSc12}
N.~Guillen and R.~W. Schwab, \emph{Aleksandrov-{B}akelman-{P}ucci type
  estimates for integro-differential equations}, Arch. Ration. Mech. Anal.
  \textbf{206} (2012), no.~1, 111--157. \MR{2968592}

\bibitem[KS79]{KrSa79}
N.~V. Krylov and M.~V. Safonov, \emph{An estimate for the probability of a
  diffusion process hitting a set of positive measure}, Dokl. Akad. Nauk SSSR
  \textbf{245} (1979), no.~1, 18--20.

\bibitem[Lan71]{Lan71}
E.~M. Landis, \emph{Uravneniya vtorogo poryadka ellipticheskogo i
  parabolicheskogo tipov}, Izdat. ``Nauka'', Moscow, 1971. \MR{0320507 (47
  \#9044)}

\bibitem[Mim13a]{Mi2}
A.~Mimica, \emph{Harnack inequality and {H}{\"{o}}lder regularity estimates for
  a {L}{\'{e}}vy process with small jumps of high intensity}, J. Theor. Probab.
  \textbf{26} (2013), 329--348.

\bibitem[Mim13b]{Mi3}
\bysame, \emph{On harmonic functions of symmetric {L}\'{e}vy processes}, Ann.
  Inst. H. Poincar\'{e} Probab. Statist. (2013), to appear.

\bibitem[RY05]{RY}
D.~Revuz and M.~Yor, \emph{Continuous martingales and {B}rownian motion},
  Springer, Berlin, 2005.

\bibitem[Sat99]{S}
K.-I. Sato, \emph{{L}\' evy processes and infinitely divisible distributions},
  Cambridge University Press, Cambridge, 1999.

\bibitem[Sil06]{Sil06}
L.~Silvestre, \emph{H\"older estimates for solutions of integro-differential
  equations like the fractional {L}aplace}, Indiana Univ. Math. J. \textbf{55}
  (2006), no.~3, 1155--1174. \MR{2244602 (2007b:45022)}

\bibitem[{\v{S}}SV06]{SSV06}
H.~{\v{S}}iki{\'c}, R.~Song, and Z.~Vondra{\v{c}}ek, \emph{Potential theory of
  geometric stable processes}, Probab. Theory Related Fields \textbf{135}
  (2006), no.~4, 547--575. \MR{2240700 (2008h:60319)}

\end{thebibliography}
\end{document}